\documentclass[a4paper]{article}
\usepackage[english]{babel}
\usepackage[utf8]{inputenc}
\usepackage{amsmath,amsthm}
\usepackage{amsfonts,amssymb}

\usepackage{tikz,tikz-cd} 

\usepackage{xifthen} 

\newcommand{\toposdefaut}{0}
\newcommand{\topos}[1][\toposdefaut]{ 
\ifthenelse{\equal{#1}{0}}{ \mathcal{T} }
{
\ifthenelse{\equal{#1}{1}}{ \mathcal{E} }{ #1 }
}
}

\newcommand{\sh}{\textsf{Sh}}

\newcommand{\C}{\mathbb{C}}

\newcommand{\Ecal}{\mathcal{E}} 
 
\newcommand{\Tcal}{\mathcal{T}}

\newcommand{\Ocal}{\mathcal{O}} 
\newcommand{\Pcal}{\mathcal{P}} 
 
\newcommand{\Scal}{\mathcal{S}}

\newcommand{\Lcal}{\mathcal{L}} 
\newcommand{\Mcal}{\mathcal{M}}

\usepackage{titlesec}
\titleformat{\subsubsection}[runin]{\normalfont}{\thesubsubsection}{0pt}{}[.]

\renewcommand{\thesubsubsection}{\arabic{section}.\arabic{subsubsection}}
\csname @addtoreset\endcsname{subsubsection}{section}

\newcommand{\block}[1]
{

\par \subsubsection{} #1

\bigskip}

\newcommand{\blockn}[1]{\par #1 \bigskip}
\newcommand{\blockp}[1]{\par #1}

\newcommand{\Th}[1]
	{
	\bigskip	
	\textbf{Theorem : }{\itshape #1}
		
	\bigskip
	}

\newcommand{\Prop}[1]
	{

	\bigskip
	
	\textbf{Proposition : }{\itshape #1}
		
	\bigskip
	
	}

\newcommand{\Lem}[1]
	{

	\bigskip
	
	\textbf{Lemma : }{\itshape #1}
		
	\bigskip
	
	}

\newcommand{\Dem}[1]{
	
	\smallskip
	
	\textbf{Proof : } \par
	 {#1} $\square$
	 
	 \bigskip
}

\begin{document}

\pagestyle{plain}
\title{A Geometric Bohr topos}
\author{Simon Henry}

\maketitle

\begin{abstract}

In this short note, we construct a variant of the Bohr topos of a $C^{*}$-algebra which takes into account the topology of the algebra in a finer way and such that this construction is stable under pullback along geometric morphisms. This generalizes a construction for finite dimensional algebras of G.Raynaud. Our idea is to construct the Bohr topos of a $C^{*}$-algebra $A$ as the sublocale of the lower power locale of the localic completion of $A$ which classifies the commutative localic sub-$C^{*}$-algebras of $A$.

\end{abstract}

\renewcommand{\thefootnote}{\fnsymbol{footnote}} 
\footnotetext{\emph{email:} simon.henry@imj-prg.fr}
\renewcommand{\thefootnote}{\arabic{footnote}} 


\tableofcontents

\section{Introduction and preliminaries}

\blockn{The term ``Bohr topos" denotes the following construction (see \cite{heunen2009topos}): one starts with a $C^{*}$-algebra $A$, we consider the poset $P_A$ of commutative sub-$C^{*}$-algebras of $A$ ordered by inclusion. One then defines the ``Bohr topos" ``$S_d(A)$"\footnote{The ``$d$" stand for ``discrete", by opposition to the construction introduced in this paper which caries a more subtle topology.} to be the topos of presheaves over $P_A^{op}$. $S_d(A)$ is then tautologically endowed with a presheaves of commutative $C^{*}$-algebras to which one can for example applies the internal Gelfand duality of \cite{banaschewski2006globalisation} to obtain a bundle $\Sigma \rightarrow S_d(A)$.}

\blockn{As explained in \cite{heunen2008principle} this construction is meant to be used in a topos theoretic approach of theoretical physics based on two principles: 

\begin{itemize}
\item The ``tovariance principle", that the law of physics should be formalised internally in a topos and invariant under geometric morphisms.

\item Bohr's doctrine of classical concept that ``all experience must ultimately be expressed in terms of classical concepts".
\end{itemize}
}

\blockn{From a purely mathematical point of view, this construction of the Bohr topos $S_d(A)$ has some defects, which have been highlighted in G.Raynaud's thesis \cite{raynaud2014fibred}, and which are in fact all related.

The first is that it ignore any form of topological structures on the set of commutative $C^{*}$-algebras. For example, in the simplest non-trivial case $A=M_2(\C)$ the space $S_d(A)$ we get is the set $\{*\} \cup \{$Orthogonal basis of $\C^{2} \}$, endowed with the topology whose open subsets are $\emptyset$ and all the subsets containing $*$. We would clearly prefer a topology compatible with the natural topology of the set of orthogonal basis of $\C^{2}$, and this keep happening with higher dimensional $C^{*}$-algebras.

The second is that the association $A \mapsto S_d(A)$ is itself not stable by pullback along geometric morphisms (this is actually a consequence of the fact that the topology on $S_d(A)$ is not the expected one), which, at least in my opinion, seems to be a relatively bad start with respect to the tovariance principle mentioned above.

The third is that (except in the finite dimensional case) $S_d(A)$ has way more points than it should. It has one point by commutative sub-$C^{*}$-algebra of $A$ whose specialization order corresponds to the inclusion order but it also have plenty of other points corresponding to non-principal prime ideals of $P_A$.
}

\blockn{In his thesis (\cite{raynaud2014fibred}) G.Raynaud proposed a new construction of the Bohr topos $S(A)$ when $A$ is a \emph{finite dimensional} $C^{*}$-algebra. His construction is based on the explicit description of the set of commutative sub-$C^{*}$-algebras of a finite dimensional $C^{*}$-algebra as a union of sub-spaces of Grasmanian manifolds and the acknowledgement that this set do comes with ``levelwise" natural topologies which can be patched together nicely by a result of C.Townsend (\cite{townsend2011patch}).}

\blockn{In order to generalize this construction and to solve these problems for a general $C^{*}$-algebra we follow a completely different approach, related to the idea of ``geometric logic" as presented for example in \cite{vickers2012coherent}: we will construct $S(A)$ as a classifying topos. Indeed this is the best way to ensure that $A \mapsto S(A)$ will be stable by pullback along geometric morphisms and we already know what points of $S(A)$ should be: commutative sub-$C^{*}$-algebras of $A$.  }

\blockn{Unfortunately, it is \emph{not} possible to define $S(A)$ as the classifying topos for commutative sub-$C^{*}$-algebras of $A$. The first reason is of course that this is not a geometric theory, but this only means that one can not deduce its existence from a general theorem, the deeper problem is that, as observed in the author's thesis (see \cite{henry2014localic}), $C^{*}$-algebras do not descend along open surjections while something classified by a topos should always have this descent property. In particular results of \cite{henry2014localic} show that such a ``Geometric Bohr  topos" should automatically have points correspondind to all the ``localic sub-$C^{*}$-algebras" of $A$ as well.}

\blockn{Let us briefly explain what is a localic $C^{*}$-algebra, and first what are frames, locales and pointfree topology. A frame is a poset which admits arbitrary supremums\footnote{hence also arbitrary infimums} and such that binary infimums distribute over arbitrary supremums, and a frame homomorphism if an order preserving map which preserve arbitrary supremums and finite infimums. In particular, a topology on a set $X$ is exactly the data of a subframe of the frame $\Pcal(X)$ of subset of $X$. The idea of locales theory (also called ``pointfree topology", or ``formal topology") is to consider that any frame is a topology independently of the fact that it is embedded in a frame of the form $\Pcal(X)$ or not. The category of locales is hence defined as the opposite of the category of frames. This category is extremely close to the category of topological spaces: there is an equivalence of categories between spatial\footnote{Those admiting sufficiently many ``points" i.e. morphisms from the terminal locale.} locales and sober\footnote{This include in particular all Hausdorff topological spaces, and any underlying topological space of a scheme} topological spaces. But the category of locales is slightly better\footnote{We refer to \cite{johnstone1983point} for an explanation of this claim.} behaved, especially in relation to topos theory, for example it is a reflexive full-subcategory of the category of toposes. An informal introduction to the theory of locales can be found in \cite{johnstone1983point}, a more complete one in \cite{picado2012frames} or in the part $C$ of \cite{sketches}. }

\blockn{A localic $C^{*}$-algebra is then a $C^{*}$-algebra whose ``underlying set" is not a set but a locale and the structure maps (addition, multiplication, involution and norm) are morphism in the category of locales. The theory of localic $C^{*}$-algebras (also called $C^{*}$-locales) has been developed in \cite{henry2014localic}, including a version of the (constructive) Gelfand duality for these algebras. Any ordinary $C^{*}$-algebra can be completed into a localic $C^{*}$-algebra (called its localic completion) and this construction induces an equivalence between the category of ordinary $C^{*}$-algebras and the category of ``weakly spatial" localic $C^{*}$-algebras (those which have a fiberwise dense set of points). Assuming the axiom of choice\footnote{Dependant choices is enough.} any localic $C^{*}$-algebra is weakly spatial, but on an arbitrary topos this is not the case. In fact localic $C^{*}$ algebras in an arbitrary topos $\Tcal$ correspond to general bundles of $C^{*}$-algebras over $\Tcal$, while ordinary $C^{*}$-algebras of $\Tcal$ correspond to bundles of $C^{*}$-algebras which have enough locale continuous sections.

We refer the reader to \cite{henry2014localic} for more details on the theory of localic $C^{*}$-algebras, as well as for some basic preliminaries on the the theory of locales, on locally positive locale and fiberwise closedness that can be of use of the understanding certain details of the present paper.}

\blockp{The main result of this paper is the following theorem:}

\block{\label{MainTh}\Th{Let $A$ be a $C^{*}$-locale. Then there is locale $S(A)$ which classifies commutative sub-$C^{*}$-locales of $A$. Moreover, the association $A \mapsto S(A)$ is compatible\footnote{In the sense that if $f$ is some geometric morphism to the base topos then $S(f^{\sharp}(A)) \simeq f^{\sharp}(S(A))$. See below for the definition of $f^{\sharp}$.} with pull-back along geometric morphisms.}

This means that for any locale (or topos) $X$, with $p$ the canonical geometric morphism from $\sh(X)$ to the base topos there is a bijection between morphisms from $X$ to $S(A)$ and commutative sub-$C^{*}$-locales of $p^{\sharp}(A)$ internally in $\sh(X)$.

In this theorem, the algebra $A$ can be unital or not, and one can consider sub-algebras containing $1$ or arbitrary sub-algebras.
}

\blockn{Section \ref{secPowerlocale} review a known construction of the theory of locales: the lower power locale. It is a generalization of hyperspace construction in classical ``point-set" topology which will be the key step in our construction of $S(A)$. In section \ref{secProof} we construct $S(A)$, hence proving our main theorem \ref{MainTh} and we give some of its basic properties. }

\blockn{In all of this paper we are working internally in an elementary topos $\Scal$ admitting a natural number object. We simply call its object sets, and its sub-object classifier is denoted $\Omega$. We do not need to assume that the base topos satisfies the law of excluded middle or the axiom of choice. In particular, the theorem stated above will be valid in this context. }

\blockn{An open sublocale $U$ of a locale $X$ (i.e. an element of the corresponding frame $\Ocal(X)$) is said to be positive if whenever one has $U = \bigcup_{i \in I} U_i$ for some family of open sublocales $(U_i)_{i \in I}$ then $\exists i \in I$. This is a positive way of saying that $U$ is non-zero. A locale $X$ is said to be locally positive if any open sublocale of $X$ can be written as a union of positive open sublocale. Assuming the law of excluded middle any locale is locally positive, but in general this is a non-trivial property: a locale is locally positive if and only its map to the terminal locale is an open map. More details and references about this can be found in the preliminaries section of \cite{henry2014localic}. }

\blockn{We also refer the reader to the preliminaries section of \cite{henry2014localic} or to \cite[C1.1 and C1.2]{sketches} for the notion of fiberwise closed sublocale (also called weakly closed sublocale). }

\blockn{A sup-lattice is a poset having arbitrary supremums and sup-lattice morphisms are the map preserving these supremums. If $f : \Tcal \rightarrow \Ecal$ is a geometric morphism between two toposes and $S$ is a sup-lattice object in $\Tcal$ (i.e. a poset object which is internally a sup-lattice) then $f_*(S)$ is also a sup-lattice and this defines a functor between sup-lattices of $\Tcal$ and sub-lattices of $\Ecal$. This functor has a left adjoint denoted $f^{\sharp}$ which can be defined in the following way : starting from a sup-lattice $S$ in $\Ecal$ one takes a presentation of it (as a sup-lattice), one pullbacks this presentation to $\Tcal$ and one construct the sup-lattice in $\Tcal$ generated by this pulled-back presentation. If $F$ is a frame in $\Ecal$ then $f^{\sharp}(F)$ will be a frame in $\Tcal$. We also denote by $f^{\sharp}$ the same functor from the category of locales of $\Ecal$ to locale of $\Tcal$ (i.e. on the opposite category). Finally, locales internal to the topos $\Ecal$ correspond to toposes endowed with a localic geometric morphism to $\Ecal$, and under this correspondence $f^{\sharp}$ corresponds to the categorical pullback.}

\blockn{We thanks S.Vickers for sugesting this problem to us as well as for mentioning the works of his student G.Raynaud. }

\section{Recall on the lower power locale}
\label{secPowerlocale}

\block{Let us start by the following observation:
\label{thP_L}

\Th{Let $\Lcal$ be any locale, then there exists a locale $P_L(\Lcal)$ which classifies locally positive fiberwise closed sublocales of $\Lcal$. Moreover the association $\Lcal \mapsto P_L(\Lcal)$ is compatible with pullback along geometric morphisms.}

$P_L(\Lcal)$ is called the lower power locale of $\Lcal$ and this theorem is proved in \cite{bunge1996constructive}, but we will need to know a little more about this result and the construction of $P_L(\Lcal)$.}

\block{The proof of the previous theorem goes through the following proposition:

\Prop{Let $\Lcal$ be any locale, then there is a bijection between fiberwise closed locally positive sublocales of $\Lcal$ and sup-lattice morphisms $\Ocal(\Lcal)$ to $\Omega$, given by the following:

\begin{itemize}

\item To any locally positive fiberwise closed sublocale $F \subset \Lcal$ one associate the map $U \mapsto ``U \wedge F \text{ is positive} "$.

\item To any sup-lattice morphisms $f: \Ocal(X) \rightarrow \Omega$ one associates the classifying locale for points $x$ of $X$ such that $(x \in U) \Rightarrow f(U)$.

\end{itemize}
}

This is also proved in \cite[section 3]{bunge1996constructive}.
}

\block{The forgetful functor from the category of frames to the category of sup-lattices has a left adjoint $\Sigma$ (the ``free frame functor"). If one defines $P_L(X)$ by:

\[\Ocal(P_L(X)) = \Sigma \Ocal(X)  \]

Then the points of $P_L(X)$ are the frame homomorphisms from $\Sigma \Ocal(X)$ to $\Omega$ which are exactly the sup-lattice maps from $\Ocal(X)$ to $\Omega$ which by the previous proposition are exactly the fiberwise closed locally positive sublocales of $X$.
}

\blockn{The final observation to conclude the proof of \ref{thP_L} is that the construction $P_L$ is compatible with pullback along geometric morphisms, but this follows immediately from its definition and the the observation that the pullback of a locale and of the underlying sup-lattice are the same functor.
}

\block{\label{PLfunct}If $f:\Lcal \rightarrow \Mcal$ is a morphism of locales, then $f^{*} : \Ocal(\Mcal) \rightarrow \Ocal(\Lcal)$ can be seens as a sup-lattice homorphism and hence defines a frame homomorphism $\Sigma f : \Sigma \Ocal(\Mcal) \rightarrow \Sigma \Ocal(\Lcal)$, i.e. a morphism of locales $P_L(f) : P_L(\Lcal) \rightarrow P_L(\Mcal)$.

Translating this definition of $P_L(f)$ in term of sublocales one obtains that:

\Lem{ $P_L(f)$ acts on generalized points by the sending any locally positive fiberwise sublocale of $\Lcal$ to the fiberwise closure of its (regular) image in $\Mcal$. }

\Dem{As the construction of $P_L(f)$ clearly commutes to pullback along geometric morphisms it is enough to check this on ordinary points. if $F$ is a fiberwise closed locally positive sublocale of $\Lcal$ corresponding to a point $p$ of $P_L(\Lcal)$, then its image by $f$ corresponds by construction to the map from $\Ocal(\Mcal)$ to $\Omega$ which send any $U$ to ``$f^{*}(U) \cap F$ is positive". This is equivalent to the fact that $U \cap f_!(F)$ is positive which is equivalent to the fact that the intersection of $U$ with the fiberwise closure of $f_!(F)$ is positive, and hence this point does correspond to the fiberwse closure of the image of $F$ by $f$. }

}

\block{ \label{PLsubloc} \Prop{Let $i :X \hookrightarrow Y$ be the inclusion of a sublocale, then $P_L(i)$ is also an inclusion of sublocale.}

\Dem{In terms of frames one has to prove that if $f : A \rightarrow B$ is a surjective frame homomorphism then $\Sigma f$ is also surjective. The image of a frame homomorphism is a subframe, indeed, because any point in the image has a canonical pre-image (the supremum of all its pre-images), a supremum (or an intersection) of points in the image is the image of the supremum (or the intersection) of the canonical pre-image. Hence the image of $\Sigma f $ is a subframe of $\Sigma B$. But it contains $B$ because it is the image of $A \subset \Sigma A$ by $\Sigma f$ and $B$ generates $\Sigma B$ as a frame hence $\Sigma f$ is surjective.
}
}

\section{The construction and first properties of $S(A)$}
\label{secProof}
\block{We will now prove the theorem \ref{MainTh}. We will directly construct $S(A)$ as a classifying space, hence the compatibility with pullback will be automatic. We fix a localic\footnote{If one start with an ordinary $C^{*}$-algebra, then $A$ will denotes its localic completion.} $C^{*}$-algebra $A$, and we will construct $S(A)$ as a sublocale of $P_L(A)$. Before that, we need two additional constructions:}

\block{For any locale $\Lcal$ there is a classifying space $P_L^{(2)}(\Lcal)$ for pairs $X \subset Y \subset \Lcal$ of fiberwise closed locally positive sublocales of $\Lcal$. The locale $P_L^{(2)}(\Lcal)$ is a sublocale of $P_L(\Lcal) \times P_L(\Lcal)$ containing the diagonal. It can be defined either as a lax-pullback (because the inclusion of closed sub-locales corresponds to the specialization order) or more directly by applying $P_L$ to the universal locally positive fibewise closed sublocale of $\Lcal$ internally in $P_L(\Lcal)$. }

\block{Also, there is a map $P_L(\Lcal) \times P_L(\Mcal) \rightarrow P_L(\Lcal \times \Mcal)$ which (on generalized points) send two fiberwise closed locally positive sublocales $X \subset \Lcal$ and $Y \subset \Mcal$ to $X \times Y \subset \Lcal \times \Mcal$. }

\block{We are now ready to prove our theorem.
\Dem{As mentioned above, we will construct $S(A)$ as a sublocale of $P_L(A)$. Indeed, any sub-$C^{*}$-locale of $A$ is in particular a fiberwise closed (because it is complete) locally positive (because it is metric) sublocale of $A$. Moreover, the difference between a commutative sub-$C^{*}$-locale and an arbitrary fiberwise closed locally positive sublocale is just a series of purely algebraic axioms that can all be translated into finite projective limits involving the structure mentioned above.

More precisely:
\begin{itemize}

\item A generalized point $X \rightarrow P_L(A)$ corresponds to a $*$-stable sublocale if and only it factor into the equalizer of the identity map and the map $P_L(*)$. The same thing in true for the stability by opposite.

\item A generalized point $X \rightarrow P_L(A)$ corresponds to a sublocale stable by $+$ if and only if it factor into the sublocale defined as the pullback of $P_L^{(2)}(A)$ by the map:

\[ P_L(A) \overset{(p,Id)}{\rightarrow} P_L(A) \times P_L(A) \]

Where $p$ is the map:

\[ P_L(A) \overset{(Id,Id)}{\rightarrow} P_L(A) \times P_L(A) \rightarrow P_L(A \times A) \overset{P_L(+)}{\rightarrow} P_L(A). \]

The same thing holds for the stability by multiplication, and one can do something similar for the stability by complex multiplication.

\item The commutativity of the sub-algebra is just slightly more involved. One can define a sublocale $C \subset A \times A$ which classifies the pair of elements commuting together. Because $A$ is metric, it is fiberwise separated and hence $C$ is fiberwise closed in $A \times A$. By \ref{PLsubloc} (and the geometric description of the functoriality of $P_L$ given in \ref{PLfunct} ) one has a sublocale $P_L(C) \subset P_L(A \times A)$ which classifies exactly the locally positive fiberwise closed sublobale of $A \times A$ included in $C$. 

And hence, a generalized point $ f :X \rightarrow P_L(A)$ corresponds to a sublocale of $A$ on which the restriction of the multiplication is commutative if and only it factor into the pullback of $P_L(C)$ by the map:

\[ P_L(A) \overset{\Delta}{\rightarrow} P_L(A) \times P_L(A) \rightarrow P_L(A \times A) \]

\item Finally, if $A$ is unital and if one wants to consider only unital sub-algebras of $A$, then there is a point $p$ of $P_L(A)$ corresponding to the sublocale $\{1\}$ (which is fiberwise closed and locally positive) of $A$, and the pullback of $P_L^{(2)}(A) \subset P_L(A) \times P_L(A)$ along $(p,id_{P_L(A)})$ classifies the fiberwise closed locally positive sublocale of $A$ containing $1$.

\end{itemize}

Hence if one defines $S(A)$ as the intersection of all these sublocales of $P_L(A)$ one obtains that a generalized point $X \rightarrow P_L(A)$ factor into $S(A)$ if and only if the corresponding sublocale of $A$ (internally in $X$) is a (unital or not) commutative sub-$C^{*}$-algebra of $A$. And hence $S(A)$ do satisfy the universal property of the theorem. The compatibility with pullback follows immediately from this universal property.
}

}

\blockn{Although this construction of $S(A)$ might seem non-explicit it actually gives a basis of the topology of $S(A)$: as $S(A)$ is a sublocale of $P_L(A)$ a basis of its topology is given by $(W_U)_{U \in \Ocal(A)}$ where $W_U$ is the sub-locale of $S(A)$ which classifies the commutative sub-$C^{*}$-locale of $A$ whose intersection with $U$ is positive. Assuming classical logic this gives a completely explicit desciption at least of the topological space of points of $S(A)$. We do not know if, assuming classical logic and the axiom of choice, $S(A)$ is spatial or not. }

\blockn{As with the previous construction of a Bohr topos, $S(A)$ is tautologically endowed with a bundle of commutative $C^{*}$-algebras, which is a sub-bundle of the constant bundle $p^{\sharp}(A)$ (where $p$ is the canonical map $S(A)\rightarrow \{* \}$). The difference with the previous construction is that now this bundle is internally a localic $C^{*}$-algebra instead of an ordinary $C^{*}$-algebra. As explained in the introduction, this just mean that this bundle of algebras might fail to have enough continuous locale sections (and hence cannot be studied through its sheaf of locale sections which corresponds to the ordinary $C^{*}$-algebra of points).

There is a theorem asserting that over a paracompact/locally paracompact basis any bundle of $C^{*}$-algebras (and more generally of Banach spaces) admit enough continuous section . It is proved in the appendix of \cite{fell1977induced} for classical bundle theory and in the author's thesis \cite[chapter 3, section 5]{henry2014thesis} for the topos theoretic version. This explain why the finite dimensional case works without using localic algebras but it is unlikely that a space like $S(A)$ do satisfies this kind of paracompactness hypothesis outside of the finite dimensional case.
}

\blockn{Thanks to our results in \cite{henry2014localic} one can still apply the constructive Gelfand duality and one do obtain a map $\Sigma_A \rightarrow S(A)$ which is internally the spectrum of this bundle of commutative localic $C^{*}$-algebras. If we are working with unital $C^{*}$-algebras then it is a proper and separated map. In the non-unital case, the gelfand duality for both ordinary and localic $C^{*}$-algebra has been covered in \cite{henry2014nonunital}, and hence one still get a map $\Sigma_A \rightarrow S(A)$ but this map will not be proper (it satisfies some fiberwise locale compactness condition instead).}

\blockn{Because of the nice ``geometric" description of $S(A)$ it is now trivial to describe what $\Sigma_A$ classifies: it classifies couple $(B,\chi)$ where $B$ is a commutative sub-$C^{*}$-locale of $A$ and $\chi$ is a non-degenrate character of $B$. Assuming classical logical, this implies that points of $\Sigma_A$ are couple $(B,\chi)$ where $B$ is a commutative sub-algebra of $A$ and  $\chi$ is a non zero character of $B$.} 

\blockn{Finally, our approach relates to that of Heunen, Landsman and Spitters in the following way: because $S(A)$ has a point of each commutative sub-$C^{*}$-locales of $A$, with the correct specialization order, there is a canonical map from $S_d(A)$ to $S(A)$. Moreover, the bundle of commutative $C^{*}$ algebras over $S_d(A)$ (and its Gelfand spectrum) studied in \cite{heunen2008principle}, \cite{heunen2009topos} are by definition the pullback of those on $S(A)$ by this map.}

\bibliography{Biblio}{}
\bibliographystyle{plain}

\end{document}